\documentclass[reqno,centertags,11pt,draft]{amsart}
\usepackage{amsmath,amsthm,amsfonts,amssymb,enumerate}
\usepackage[bookmarksopen=true,final]{hyperref}
\usepackage{multirow,yhmath}
\usepackage{xcolor}
\usepackage{bm}
\usepackage[nobysame]{amsrefs}

\newcommand{\bbC}{{\mathbb{C}}}
\newcommand{\bbD}{{\mathbb{D}}}

\newcommand{\bbN}{{\mathbb{N}}}

\newcommand{\bbR}{{\mathbb{R}}}

\newcommand{\bbT}{{\partial\mathbb{D}}}

\newcommand{\bbZ}{{\mathbb{Z}}}

\newcommand{\calM}{{\mathcal M}}
\newcommand{\calN}{{\mathcal N}}

\newcommand{\la}{\left\langle}
\newcommand{\ra}{\right\rangle}

\newcommand{\supp}{\text{\rm{supp}}}

\newcommand{\beq}{\begin{equation}}
\newcommand{\eeq}{\end{equation}}
\newcommand{\ba}{\begin{align*}}
\newcommand{\ea}{\end{align*}}

\newcommand{\abs}[1]{\lvert#1\rvert}

\DeclareFontFamily{U}{mathx}{}
\DeclareFontShape{U}{mathx}{m}{n}{<-> mathx10}{}
\DeclareSymbolFont{mathx}{U}{mathx}{m}{n}
\DeclareMathAccent{\widehat}{0}{mathx}{"70}
\DeclareMathAccent{\widecheck}{0}{mathx}{"71}

\DeclareMathOperator{\re}{Re}

\allowdisplaybreaks
\numberwithin{equation}{section}

\usepackage[normalem]{ulem}

\newtheorem{theorem}{Theorem}
\newtheorem{lemma}[theorem]{Lemma}

\newtheorem{definition}[theorem]{Definition}

\newtheorem{remark}[theorem]{Remark}

\begin{document}
\title{Nikishin systems on the unit circle}
\author{Rostyslav Kozhan}
\date{\today}
\begin{abstract}
	We introduce Nikishin system of $r$ probability measures on the unit circle. Using a generalized Andreief identity, we show that such systems satisfy the AT property and therefore normality, introduced in~\cite{KVMLOPUC}, for any multi-index 
	$\bm{n}\in\bbN^r$ with same-parity components satisfying $n_1 \ge n_2 \ge\ldots\ge n_r$. 
	In the case of $r=2$, we demonstrate that the same property holds without requiring $n_1 \ge n_2 \ge\ldots\ge n_r$.
	
	The analogous simple proof works for Nikishin systems on the real line for indices satisfying $n_j\ge \max\{n_{j+1},\ldots,n_r\}-1$, $j=1,\ldots,r-1$. This is related to the proof by Cousseement and Van Assche~\cite{CouVA} for $r=2$.
\end{abstract}
\maketitle
\section{Introduction}

Given a probability measure $\mu$ on the real line $\bbR$ or on the complex unit circle $\bbT$, one can associate a sequence of monic polynomials $\{P_n\}_{n=0}^\infty$ with $\deg P_n = n$ 
that are pairwise orthogonal in $L^2(\mu)$:
\begin{equation*}
	\int_\bbR P_{n}(x) {P_m(x)} d\mu(x) = 0, \qquad \mbox{if } n\ne m,
\end{equation*}
for measures on the real line and
\begin{equation*}
	\int_\bbT \Phi_{n}(z) \overline{\Phi_m(z)} d\mu(z) = 0, \qquad \mbox{if } n\ne m,
\end{equation*}
for measures on the unit circle.
Behavior of these polynomials, their zeros, and their recurrence coefficients 
carry a lot of useful information about the original
measure $\mu$. As a result, theory of orthogonal polynomials is an important and versatile tool in many branches of mathematics as well as in physics and engineering. We refer the reader to~\cite{SimonL2,Chihara} for the theory of orthogonal polynomials on the real line and to~\cite{OPUC1,OPUC2} for the the theory on the unit circle. 

If we require simultaneous orthogonality conditions with respect to several orthogonality measures $\mu_1,\ldots,\mu_r$, then we arrive at the notion of \textit{multiple orthogonal polynomials}. The theory of such polynomials is very well-developed on the real line and finds applications in approximation theory, analytic number theory, random matrix theory, and integrable probability. 
See \cite{Applications} for a quick introduction, \cite{Ismail,bookNS} for a comprehensive treatment, and~\cite{Aptekarev,VA-HP,VA99} for reviews.

One of the fundamental notions in this theory is that of \textit{normality}. We say that $(n_1,\ldots,n_r)\in\bbN^r$ is normal with respect to $(\mu_1,\ldots,\mu_r)$, if there exists a \textit{unique} monic univariate polynomial $P_{\bm{n}}$ of exact degree $n_1+\ldots+n_r$ satisfying 
simultaneous orthogonality conditions
	\begin{equation}\label{eq:moprlIIIntro}
	\int P_{\bm{n}}(x)x^p d\mu_j(x) = 0,\qquad p = 0,1,\dots,n_j-1 ,\qquad j = 1,\dots,r.
\end{equation} 
We say that a system is \textit{perfect}, if every $\bbN^r$ index is normal. 


In~\cite{Nik80}, Nikishin introduced a class of systems  $(\mu_1,\ldots,\mu_r)$ on $\bbR$  that are now called the Nikishin systems (see Definition~\ref{def:Nikishin}). These systems are natural from the point of view of approximation theory
and were shown to be perfect. Many authors contributed to the proof of this statement. 
Initially Nikishin ~\cite{Nik80} showed that every \textit{step-line} index (i.e., index of the form $(j+1,\ldots,j+1,j,\ldots,j)$) is normal. 
Driver and Stahl in~\cite{DriSta94} 
extended this to any index $\bm{n}\in\bbN^r$ that satisfies
\begin{equation}\label{eq:simpleNIntro}
	n_j\ge \max\{n_{j+1},\ldots,n_r\}-1, \qquad \mbox{for any } j=1,\ldots,r-1.
\end{equation}
The set of indices that were normal was further extended  in~\cite{BBFL,FidIllLop04}.
It was shown that \textit{every} index is normal if $r=2$ in~\cite{DriSta94}, if $r=3$ in~\cite{FidLop02}. 
Finally, in~\cite{FidLop11}, Fidalgo and L\'{o}pez showed this for any $r$, and further improved this in~\cite{FidLop11b} by allowing the supports to be touching and/or unbounded. 

The first result of this paper is a simplified proof of normality of indices~\eqref{eq:simpleNIntro} for any Nikishin system which allows the supports to be touching and/or unbounded (see Section~\ref{ss:Nikishin}). The proof was inspired by the proof of Coussement and Van Assche~\cite{CouVA} for $m=2$ (see also~\cite{Kui}). We streamline with the help of a generalization of the Andreief identity (see Lemma~\ref{lem:det}\textit{(i)}) and generalize it to any $m$ using an inductive argument.  

The main aim of this paper is to introduce Nikishin systems \textit{on the unit circle}. Multiple orthogonal polynomials on the unit circle were introduced by M\'{i}nguez and Van Assche in~\cite{MOPUC1}, followed by the papers~\cite{MOPUC2,KVSzego}. In~\cite{KVMLOPUC}, the authors introduces \textit{Laurent} multiple orthogonality (see Definition~\ref{def:MLOPUC}) which allowed to naturally introduce the notion of Angelesco and AT systems and to prove their perfectness. They showed that this setting is natural as it relates to a two-point Hermite--Pad\'{e} approximation of Carath\'{e}odory functions.

In this paper, we introduce the analogue of the Nikishin systems on the unit circle. Using the generalized Andreief identity and an inductive argument as on $\bbR$, we show that any index $(n_1,\ldots,n_r)\in\bbN^r$ with  same-parity $n_j$'s and with
\begin{equation}\label{eq:simpleNopucIntro}
	n_1\ge n_2 \ge \ldots \ge n_r,
\end{equation}
is normal (compare with~\eqref{eq:simpleNIntro}), see Section~\ref{ss:NikishinOPUC}.

In Section~\ref{ss:NikishinOPUCr=2}, we consider the case $r=2$ and show the normality of \textit{every} index with same-parity components without the restriction~\eqref{eq:simpleNopucIntro}. Whether the same holds for $r>2$ remains an interesting open question.

\section{Preliminaries}
\subsection{Notation}
\hfill\\

Let $\bbR$ be the set of reals, $\bbT$ be the complex unit circle, and $\bbN:=\{0,1,2,\ldots\}$ be the set of nonnegative integers. For a multi-index $\bm{n}\in\bbN^r$ we define $|\bm{n}|:=n_1+\ldots+n_r$.


Given a set $\Delta\subseteq\bbR$ or $\Delta\subseteq\bbT$ , let  $\calM(\Delta)$ be the set of all finite Borel measures $\mu$ of constant sign (so either positive or negative), with all the moments finite, and with an infinite support  $\supp\,\mu\subseteq \Delta$. 


For any set $S\subseteq\bbR$ or $S\subseteq \bbT$, we denote  $\mathring{S}$ to be the interior of $S$ in the topology of $\bbR$ or $\bbT$, respectively.

\subsection{Basics of orthogonal polynomials on the real line (OPRL)}
\hfill\\

We remind some basics of the theory of orthogonal polynomials on the real line. For more details, see, e.g.,~\cite{SimonL2}. Given a probability measure $\mu\in\mathcal{M}(\bbR)$, 
one can form a sequence of monic orthogonal polynomials $\{P_n(x)\}_{n=0}^\infty$ satisfying
\begin{equation}\label{eq:OPRL}
\int P_n(x) x^j d\mu(x) = 0, \qquad 0\le j \le n-1,
\end{equation}
and $\deg P_n = n$. 
They satisfy the three-term recurrence relation
\begin{equation}\label{eq:threeterm}
xP_n(x) = P_{n+1}(x) + b_n P_n(x) + a^2_{n-1} P_{n-1}(x)
\end{equation}
with the recurrence coefficients $a_n,b_n\in\bbR$, which are called the Jacobi coefficients. 

To each $\mu$ we can associate the sequence of moments $c_j = \int x^j d\mu(x)$, $j\ge 0$, and the Hankel matrices $H_{n}:=(c_{i+j})_{i,j=0}^{n-1}$ which are positive-definite since
\begin{equation}\label{eq:HankelInv}
\bm{v}^* H_n \bm{v} = \int \Big(\sum_{j=0}^{n-1} v_j x^j \Big)^2 d\mu(x) >0.
\end{equation}
The converse holds true: if all Hankel determinants $H_{n}=(c_{i+j})_{i,j=0}^{n-1}$ are positive definite then $\{c_j\}_{j=0}^\infty$ is the sequence of moments for some positive measure $\mu$. Invertibility of the \textit{generalized} Hankel determinants~\eqref{eq:D} is a non-trivial and important question at the heart of the theory of multiple orthogonal polynomials. It has been studied by numerous authors, and we revisit it briefly in Section~\ref{ss:ATNikishin} for Nikishin systems.


Let us also define the $m$-function of $\mu$ which is an analytic function on $\bbC\setminus\supp\,\mu$ defined by
\begin{equation}\label{eq:m}
	m_\mu(z) := \int_{\bbR} \frac{d\mu(x)}{x-z}.
\end{equation}
It satisfies the so-called stripping relation, see e.g.,~\cite[Thm 3.2.4]{SimonL2}, which takes the form
\begin{equation}\label{eq:strippingM}
	\frac{1}{m_{\mu}(z)} = b_1-z- a_1^2 m_{\mu^{(1)}}(z),
\end{equation}
where $a_1,b_1\in \bbR$ are the Jacobi coefficients from~\eqref{eq:threeterm}, and
\begin{equation}
	m_{\mu^{(1)}}(z) = \int_{\bbR} \frac{d\mu^{(1)}(x)}{x-z} 
\end{equation}
is the $m$-function of another probability measure $\mu^{(1)}$ which is the spectral measure of the Jacobi operator of $\mu$ but with the first row and column stripped. Using analyticity, it is easy to show that if $\supp\,\mu$ belongs to some interval $\Delta$ then $\supp\,\mu^{(1)}$ also belongs to $\Delta$. 

\subsection{Basics of orthogonal polynomials on the unit circle (OPUC)}\label{ss:OPUC}
\hfill\\

The theory of orthogonal polynomials on the unit circle mirrors, to some extent, the theory on the real line. Here, we introduce a few key concepts of this theory; for proofs and additional results, interested readers may refer to, for instance,~\cite{OPUC1,OPUC2}. Given a probability measure $\mu\in\mathcal{M}(\bbT)$, 
one can form a sequence of monic orthogonal polynomials $\{\Phi_n(x)\}_{n=0}^\infty$ satisfying
\begin{equation}\label{eq:OPUC}
\int_{\bbT} \Phi_n(z) z^{-j} d\mu(z) = 0, \qquad 0\le j \le n-1,
\end{equation}
and $\deg \Phi_n = n$. 

To each $\mu$ we can associate the sequence of moments $c_j = \int z^j d\mu(z)$, $j\in\bbZ$, and the (self-adjoint) Toeplitz matrices $T_{n}:=(c_{k-j})_{j,k=0}^{n-1}$ which are positive-definite since
\begin{equation}\label{eq:ToeplitzInv}
\bm{v}^* T_n \bm{v} = \int_{\bbT} \Big| \sum_{j=0}^{n-1} v_j z^j \Big|^2 d\mu(z) >0.
\end{equation} 
For the converse, if all Toeplitz determinants $T_{n}=(c_{k-j})_{j,k=0}^{n-1}$ are positive definite then $\{c_j\}_{j\in\bbZ}$ is the sequence of moments for some positive measure $\mu$ on $\bbT$. Various \textit{generalizations}  of these Toeplitz matrices (see~\eqref{eq:DOPUC}, or~\cite[Eq. (1)]{MOPUC1}) relate to the uniqueness of multiple orthogonal polynomials in the multiple orthogonality setting. For Angelesco and AT systems this was addressed in~\cite{KVMLOPUC}, and in Section~\ref{ss:NikishinOPUC} we study it for study it for Nikishin systems.



The Carath\'{e}odory function of $\mu$ is defined as the analytic function on $\bbT\setminus\supp\,\mu$ given by
\begin{equation}\label{eq:caratheodory function}
	F_\mu(z) = \int_\bbT \frac{w+z}{w-z}d\mu(w). 
\end{equation}
It is an analytic map from $\bbD$ to $\{z\in\bbC:\re z>0\}$. 
The important property of $F_\mu$ is $\overline{F_\mu(1/\bar{z})} = -F_\mu(z)$ which implies
\begin{equation}\label{eq:CaraReal}
	iF_\mu(z)\in \bbR, \qquad z\in\bbT\setminus\supp\,\mu.
\end{equation}

It is not hard to see that if $F_\mu$ is a Carath\'{e}odory function then so is $1/F_\mu(z)$. Therefore there exists a probability measure $\mu^{(-1)}$ on $\bbT$ such that
\begin{align}
	\label{eq:CaraInverse}
	1/F_{\mu}(z) & = F_{\mu^{(-1)}}(z) 
	\\
	F_{\mu^{(-1)}}(z) &= \int_\bbT \frac{w+z}{w-z}d\mu^{(-1)}(w). 
\end{align}
The $(-1)$ index reflects that the recurrence coefficients of $\mu$ and of $\mu^{(-1)}$ are negative of each other. In fact, $\mu^{(-1)}$ is an element of a one-parameter  family of probability measures $\mu^{(\lambda)}$, $\lambda\in\bbT$ on $\bbT$, called the Aleksandrov measures, see~\cite[Sect. 1.3.9, 3.2]{OPUC1} for their properties. In particular,~\cite[Thm 3.2.16]{OPUC1} says that if $\supp\,\mu$ is a subset of some arc $\Gamma\subset\bbT$, then $\supp\,\mu^{(-1)}$ also belongs to $\Gamma$ with the possible exception that $\mu^{(-1)}$ can have one point mass on $\bbC\setminus\Gamma$ at the unique location of $\bbC\setminus\Gamma$ where $F_\mu(z)=0$.

We will also need the $m$-function, which is related to $F_\mu$:
\begin{equation}\label{eq:FandM}
	m_{\mu}(z): = \int \frac{d\mu(w)}{z-w} =\frac{1-F_\mu(z)}{2z}, \quad z\notin \supp\,\mu.
\end{equation}

\subsection{Some determinants}
\hfill\\

In what follows we make use of some of the determinantal identities which we collect here for the reader's convenience.

\begin{lemma}\label{lem:det}
    \begin{enumerate}[(i)]
    	\item (The generalized Andreief identity)
    	 Let $f_j,g_j \in L^2(\mu)$ for some probability measure on a measure space $(X,\Sigma,\mu)$.  Then for any $N\ge M\ge 1$  and any $(N-M)\times N$ matrix $\bm{A}$:
    	 \begin{multline}\label{eq:Andreief}
    	 	\det\left(
    	 	\begin{array}{@{}c@{}}
    	 		\bm{A} \\
    	 		\hline
    	 		\begin{pmatrix}
    	 			\int_X f_j(t) g_k(t)\,d\mu(t)
    	 		\end{pmatrix}_{1\le j \le M, 1\le k\le N} 
    	 	\end{array}
    	 	\right) =
    	 	\\
    	 	\tfrac{1}{M!}
    	 	\int_{X^M}
    	 	\det\left(
    	 	\begin{array}{@{}c@{}}
    	 		\bm{A}\\
    	 		\hline 
    	 		\begin{pmatrix}
    	 			g_k(y_j)
    	 		\end{pmatrix}_{1\le j \le M, 1\le k\le N}
    	 	\end{array}
    	 	\right)
    	 	\det\left(f_l(y_j)\right)_{1\le l,j \le M} \,d^M\mu(\bm{y}),
    	 \end{multline}
    	where $d^M\mu(\bm{y}):=d\mu(y_1)\ldots d\mu(y_M)$. The matrix on the left-hand side refers to the matrix whose upper $(N-M)\times N$ block is $\bm{A}$ and whose lower
    	$M\times N$ block is 
    	$\begin{pmatrix}
    		\int_X f_j(t) g_k(t)\,d\mu(t)
    	\end{pmatrix}_{1\le j \le M, 1\le k\le N}$.



        \item (The Cauchy--Vandermonde determinant) Let $\Delta_1$ and $\Delta_2$ be two closed intervals with $\Delta_1\cap\Delta_2 = \varnothing$. If $t_j\in\Delta_1$ for each $1\le j \le n_1+n_2$ and $z_j\in \Delta_2$ for each $1\le j\le n_2$, then
        \begin{multline*}
        \det
		\begin{pmatrix}
			1&1& \ldots & 1\\
			t_1 & t_2 &\ldots & t_{|\bm{n}|} \\
			\vdots & \vdots & \ddots& \vdots \\
			t_1^{n_1-1} & t_2^{n_1-1} & \ldots & t_{|\bm{n}|}^{n_1-1}
			\\
			\frac{1}{t_1-z_1}  & \frac{1}{t_2-z_1} & \ldots &   \frac{1}{t_{|\bm{n}|}-z_1}   \\
			\frac{1}{t_1-z_2}  &\frac{1}{t_2-z_2}& \ldots &   \frac{1}{t_{|\bm{n}|}-z_2} \\ 
			\vdots & \vdots& \ddots& \vdots \\
			\frac{1}{t_1-z_{n_2}}  &\frac{1}{t_2-z_{n_2}}  & \ldots &  \frac{1}{t_{|\bm{n}|}-z_{n_2}} 
		\end{pmatrix}
  \\ = (-1)^{n_1 n_2+\frac{n_2(n_2-1)}2} \prod_{1\le j<k \le |\bm{n}|} (t_k - t_j) \prod_{1\le j<k \le n_2} (z_k-z_j)  \prod_{j,k}(t_k-z_j)^{-1}.
     \end{multline*}
%
%
%
    \end{enumerate}
\end{lemma}
\begin{proof}
	{\textit{(i)}} This is proved in~\cite{KVMLOPUC}. 
	
{\textit{(ii)}}  This determinant is the so-called Cauchy--Vandermonde determinant, as it generalizes both the Vandermonde and the Cauchy determinants. 
To get a proof, subtract the first column from each of the other columns, and then extract a factor of $(t_j-t_1)$ from each column and a factor of $(t_1-z_j)^{-1}$ from the last $n_2$ rows, then proceed by induction. See~\cite{GMM} for even more general determinants of this type.
%
\end{proof}

\section{Multiple orthogonal polynomials on the real line (MOPRL)}
\subsection{Main definitions}
\hfill\\

Let $r\ge 1$ and consider a system of measures  on the real line $\bm{\mu} = (\mu_1,\dots,\mu_r) \in \mathcal{M}(\bbR)^r$. 
\begin{definition}
	We say that a multi-index $\bm{n}\in\bbN^r$ is normal with respect to $\bm{\mu}$ if there exists a unique monic polynomial $P_{\bm{n}}(x)$ with $\deg{P_{\bm{n}}} =\abs{\bm{n}}$, and
	\begin{equation}\label{eq:moprlII}
		\int P_{\bm{n}}(x)x^p d\mu_j(x) = 0,\qquad p = 0,1,\dots,n_j-1 ,\qquad j = 1,\dots,r.
	\end{equation} 
	We call $P_{\bm{n}}$ the type II monic multiple orthogonal polynomial.
\end{definition}

It is easy to see that normality of $\bm{n}\ne \bm{0}$ is equivalent to the invertibility of the $|\bm{n}|\times |\bm{n}|$ matrix
\begin{equation}
	\label{eq:D}
	\bm{H}_{\bm{n}} = 
	 \begin{pmatrix}
		H_{n_1,|\bm{n}|}^{(1)} \\
		\vdots \\
		H_{n_r,|\bm{n}|}^{(r)}
	\end{pmatrix},
\end{equation}
where $H_{n_j,|\bm{n}|}^{(j)}$ is the $n_j \times |\bm{n}|$ matrix
\begin{equation}
	H_{n_j,|\bm{n}|}^{(j)} = 
	\begin{pmatrix}
				c_0^{(j)} & c_1^{(j)} & \cdots & c_{\abs{\bm{n}}-1}^{(j)} \\
				c_1^{(j)} & c_2^{(j)} & \cdots &  c_{\abs{\bm{n}}}^{(j)} \\
				\vdots & \vdots & \ddots & \vdots \\
				c_{n_j-1}^{(j)} & c_{n_j}^{(j)} & \cdots & c_{\abs{\bm{n}}+n_j-2}^{(j)} \\
			\end{pmatrix},
\end{equation} 
with $c^{(j)}_k := \int x^k d\mu_j(x)$.

A system $\bm{\mu}=(\mu_1,\dots,\mu_r)$ in $\mathcal{M}^r$ is  perfect if every $\bbN^r$-index is normal.


\subsection{AT systems on the real line}\label{ss:ATNikishin}
\hfill\\

Let us introduce the notation
\begin{equation}\label{eq:ChebDet}
	U_{\bm{x}}(f_1 , f_2, \ldots, f_n) :=\det
	\begin{pmatrix}
		f_1(x_1) & \ldots & f_1(x_n) \\
		\vdots & \ddots& \vdots \\
		f_n(x_1) & \ldots & f_n(x_n) \\
	\end{pmatrix},
\end{equation}
where $\bm{x}:=(x_1,\ldots,x_n)$.

Let  $\Delta\subseteq\bbR$  be an interval. A collection of continuous on $\Delta$ functions $\{f_1,\ldots,f_n\}$ is called a Chebyshev system 
on $\Delta$ if 
\begin{equation}\label{eq:ChebPty}
U_{\bm{x}}(f_1 , f_2, \ldots, f_n) 
	\ne 0
\end{equation}
for every choice of distinct points $x_1,\ldots,x_n\in \Delta$. 
The condition~\eqref{eq:ChebPty} is invariant under arbitrary permutations of $x_j$'s, so we henceforth always assume that $x_j$'s are ordered according to
$x_1<x_2<\ldots <x_n$.

Assuming the functions $f_j$'s are real, continuity implies that~\eqref{eq:ChebPty} is equivalent to $U_{\bm{x}}$ having a constant sign over the chamber $x_1<x_2<\ldots <x_n$ of $\Delta^n$. 


\begin{definition}\label{def:AT}
	Let $\mu\in\mathcal{M}(\Delta)$ with infinite support. 
	For each $j=1,\ldots,r$, let $\mu_j\in\mathcal{M}(\Delta)$ be an absolutely continuous measure with respect to $\mu$, $d\mu_j(x)=w_j(x)d\mu(x)$.
	
	We say that $(\mu_1,\ldots,\mu_r)$ is an AT system on $\Delta$  at 
	$\bm{n}\in\bbN^r$ if
\begin{equation}\label{eq:AT}
	\{w_1, xw_1,\ldots,x^{n_1-1}w_1,w_2,xw_2,\ldots,x^{n_2-1}w_2,\ldots,w_r,xw_r,\ldots,x^{n_r-1}w_r\}
\end{equation}
is Chebyshev on $\Delta$. 
\end{definition}

\subsection{Nikishin systems on the real line: indices with non-increasing components are normal}\label{ss:Nikishin}
\hfill\\

Recall that the $m$-function of a measure is defined by~\eqref{eq:m}. 
Given two intervals $\Delta_1$ and $\Delta_2$ with $\mathring{\Delta}_1\cap\mathring{\Delta}_2=\varnothing$, and two measures $\sigma_1\in\calM(\Delta_1)$,  $\sigma_2\in\calM(\Delta_2)$, denote $\la  \sigma_1,\sigma_2 \ra$ to be the measure in $\calM(\Delta_1)$ given by
\begin{equation}\label{eq:brackets}
	d \la  \sigma_1,\sigma_2 \ra (x) := m_{\sigma_2}(x) d\sigma_1(x). 
\end{equation}
Note that $m_{\sigma_2}(x)>0$ for all $x\in\Delta_1$ if $\Delta_1<\Delta_2$ and $m_{\sigma_2}(x)<0$ for all $x\in\Delta_1$  if $\Delta_1>\Delta_2$. Therefore the right-hand side of~\eqref{eq:brackets} always defines a finite sign-definite measure on $\Delta_1$ assuming that $\Delta_1$ and $\Delta_2$ are disjoint.
If $\Delta_1$ and $\Delta_2$ share a common endpoint however, then by writing $\la  \sigma_1,\sigma_2 \ra$, we implicitly {\it assume} that ~\eqref{eq:brackets} defines a finite measure, as this is no longer automatic and depends on the behavior of $\sigma_1$ and $\sigma_2$ at that point. 

\begin{definition}\label{def:Nikishin}
We say that $\bm{\mu}=(\mu_1,\ldots,\mu_r)\in\calM(\bbR)^r$ forms a Nikishin system generated by $(\sigma_1,\ldots,\sigma_r)$ (we then write $\bm{\mu}=\calN(\sigma_1,\ldots,\sigma_r)$), if there is
a collection of intervals $\Delta_j$, $j=1,\ldots,r$ such that
\begin{equation}
	\mathring{\Delta}_j \cap \mathring{\Delta}_{j+1} = \varnothing, \qquad j=1,\ldots,r-1,
\end{equation}
and measures $\sigma_j \in\calM(\Delta_j)$, $j=1,\ldots,r$, so that
\begin{equation}
	\mu_1 = \sigma_1,  \mu_{2} = \la \sigma_1, \sigma_{2}\ra, \mu_{3} = \la \sigma_1, \la \sigma_{2},\sigma_3\ra \ra, \ldots, \mu_r = \la \sigma_1, \la \sigma_{2},\la \sigma_3,\ldots,\sigma_r\ra\ra \ra.
\end{equation}
\end{definition}

It will be convenient to 
define
\begin{equation}\label{eq:mtk}
	\begin{aligned}
		m_{t^j\mu} (x) &:=  \int \frac{t^j d\mu(t)}{t-x}, 
		\\
		m_{\la t^j\sigma_1, \sigma_{2}\ra} (x) &:=  \int_{\Delta_1} \frac{t^j m_{\sigma_2}(t) d\sigma_1(t)}{t-x} ,
		\end{aligned}
\end{equation}
for any $j\ge 0$. In particular, if $j=0$ then this is just the 
$m$-functions~\eqref{eq:m} for $\mu$ and for $\la  \sigma_1,\sigma_2 \ra$ from~\eqref{eq:brackets}, respectively.

It is well-known~\cite{FidLop11} that Nikishin systems satisfy AT property at every index $\bm{n}\in\bbN^r$, see more references in the Introduction. We provide a streamlined proof here that takes inspiration from~\cite{CouVA} incorporating a generalized Andreief identity (see Lemma~\ref{lem:det}\textit{(i)}). This combines into a short, non-technical, and transparent proof.  

\begin{theorem}\label{thm:Nikishin}
	Let $\bm{\mu}=(\mu_1,\ldots,\mu_r) = \calN(\sigma_1,\ldots,\sigma_r)$ be a Nikishin system. Suppose
\begin{equation}\label{eq:simpleN}
	n_j\ge \max\{n_{j+1},\ldots,n_r\}-1, \qquad \mbox{for any } j=1,\ldots,r-1.
\end{equation}
	Then $\bm{\mu}$ is AT on $\Delta_1$ at $\bm{n}$. 
\end{theorem}
\begin{proof}
	We prove the statement using an inductive argument. For $r=1$ any index is normal, since all the determinants $H_n$~\eqref{eq:D} are Hankel and invertible by~\eqref{eq:HankelInv} (recall that $\mu$ is infinitely supported). 	 
	Suppose the statement is proved for Nikishin systems of $r-1$ measures.	
	Let $\bm{\mu}$ be a Nikishin system with $r$ measures. Each $\mu_j$ is absolutely continuous with respect to $\sigma_1$, with the Radon--Nikodym derivatives $\frac{d\mu_j(x)}{d\sigma_1(x)}$ which are 
	$1$, $m_{\sigma_2}(x)$, $m_{\la\sigma_2,\sigma_3\ra}(x)$, $\ldots$, 
	$m_{\la\sigma_2,\ldots,\sigma_r\ra }(x)$, respectively.
	
	Let $\bm{n}$ satisfies~\eqref{eq:simpleN}. To show that AT property at $\bm{n}$ we need to show nonvanishing of 
\begin{multline}\label{eq:U}
	U_{\bm{x}}\big(
	1, x,\ldots,x^{n_1-1},
	m_{\sigma_2},xm_{\sigma_2},\ldots,	x^{n_2-1}m_{\sigma_2},
	m_{\la\sigma_2,\sigma_3\ra},xm_{\la\sigma_2,\sigma_3\ra},
	\ldots,
	\\
	x^{n_3-1}m_{\la\sigma_2,\sigma_3\ra}
	,\ldots,
	m_{\la\sigma_2,\ldots,\sigma_r\ra },xm_{\la\sigma_2,\ldots,\sigma_r\ra },
	\ldots,
	x^{n_r-1} m_{\la\sigma_2,\ldots,\sigma_r\ra }\big),
\end{multline}
	with an arbitrary (ordered) $\bm{x}\in\Delta_1^{|\bm{n}|}$.

	Now use $\frac{x^j}{t_2-x} = -\sum_{s=0}^{j-1} x^s t_2^{j-1-s} + \frac{t_2^j}{t_2-x}$ to obtain 
	\begin{align*}
		x^j m_{\la\sigma_2,\ldots,\sigma_d\ra }(x) 
		&= x^j \int_{\Delta_2} \frac{m_{\la\sigma_3,\ldots,\sigma_d\ra }(t_2)  }{t_2-x} d\sigma_2(t_2)
		\\
		&= -\sum_{s=0}^{j-1} v_{j,d,s} x^s 
		+ \int_{\Delta_2} \frac{t_2^j m_{\la\sigma_3,\ldots,\sigma_d\ra }(t_2)  }{t_2-x}d\sigma_2(t_2) 
		\\
		&= - \sum_{s=0}^{j-1} v_{j,d,s} x^s 
		+  m_{\la t^j \sigma_2,\ldots,\sigma_d\ra }(x) ,
	\end{align*}
	where 
	$v_{j,d,s} = \int_{\Delta_2} t_2^{j-1-s} d\la\sigma_2,\ldots,\sigma_d\ra(t_2)$ 
	are constants (recall the notation~\eqref{eq:mtk}).
	Since $n_1 \ge n_d-1$ for any $d$, we can perform elementary row operations in the determinant~\eqref{eq:U} to 	
	reduce it to
\begin{multline}\label{eq:U2}
	U_{\bm{x}}\big(
	1, x,\ldots,x^{n_1-1},
	m_{\sigma_2},m_{t\sigma_2},\ldots,	m_{t^{(n_2-1)}\sigma_2},
	m_{\la\sigma_2,\sigma_3\ra},m_{\la t\sigma_2,\sigma_3\ra},
	\ldots,
	\\
	m_{\la t^{(n_3-1)}\sigma_2,\sigma_3\ra}
	,\ldots,
	m_{\la\sigma_2,\ldots,\sigma_r\ra },m_{\la t\sigma_2,\ldots,\sigma_r\ra },
	\ldots,
	 m_{\la t^{(n_r-1)}\sigma_2,\ldots,\sigma_r\ra }\big).
\end{multline}
	In summary, we were able to replace factors of $x^j$ in~\eqref{eq:U} with the factor of $t_2^j$ in~\eqref{eq:U2}, which is very useful as this allows us to apply the generalized Andreief identity as follows. We treat~\eqref{eq:U2} as the left-hand side of~\eqref{eq:Andreief} with $N=|\bm{n}|$, $M=|\bm{n}|-n_1$, $N-M = n_1$, $\bm{A} = (x_k^{j-1})_{1\le j \le n_1-1,1\le k \le |\bm{n}|}$, $\mu = \sigma_2$, $g_k(t_2) = \frac{1}{t_2-x_k}$. Then the right-hand side of~\eqref{eq:Andreief} becomes
\begin{multline}\label{eq:inductive}
	\tfrac{1}{(|\bm{n}|-n_1)!}\int_{\Delta_2^{|\bm{n}|-n_1}}
	U_{\bm{x}}\Big(1,x,\ldots,x^{n_1-1},\tfrac{1}{y_1-x},\ldots,\tfrac{1}{y_{|\bm{n}|-n_1}-x}
	\Big)
	\\
	\times 
	U_{\bm{y}}\Big( 1, x,\ldots,x^{n_2-1},m_{\sigma_3},x m_{\sigma_3},\ldots,x^{n_3-1} m_{\sigma_3}, \ldots,
	\\
	m_{\la\sigma_3,\ldots,\sigma_r\ra },xm_{\la\sigma_3,\ldots,\sigma_r\ra },
	\ldots,
	x^{n_r-1} m_{\la\sigma_3,\ldots,\sigma_r\ra } \Big)
	d^{{|\bm{n}|-n_1}}\sigma_2(\bm{y}).
\end{multline}
Observe that the integrand in the last expression is invariant under permutation of $y_j$'s. Therefore we can restrict the domain of integration $\Delta_2^{|\bm{n}|-n_1}$ to the ordered chamber $y_1<y_2<\ldots<y_{|\bm{n}|-n_1}$ by simply removing the factor $\tfrac{1}{(|\bm{n}|-n_1)!}$. 

By the induction hypothesis, the system $\mathcal{N}(\sigma_2,\sigma_3,\ldots,\sigma_r)$ is AT at the location $(n_2,n_3,\ldots,n_3)$ (note that this truncated index also satisfies~\eqref{eq:simpleN}), which means that the second factor in the integrand in~\eqref{eq:inductive} has a constant sign on our ordered $\bm{y}$-chamber. The first factor also has a constant sign by Lemma~\ref{lem:det}\textit{(ii)}. Note that if $\Delta_1$ and $\Delta_2$ have one common point, then the first factor in~\eqref{eq:inductive} may in fact be zero at points where one of the $y$'s coincide with one of the $x$'s. However the integral is nevertheless nonzero since there is a choice of distinct $y_j$'s in $\Delta_2$ that do not belong to $\Delta_1$ (use $|\supp\,\sigma_2|=\infty$ and continuity of the integrand). 
This completes the induction step.
\end{proof}

\subsection{Nikishin systems on the real line, $r=2$ case: all indices are normal}
\hfill\\

There is a well-known elegant trick for the case of Nikishin systems with $r=2$ measures, that allows to show the AT property for any index without requiring~\eqref{eq:simpleN}, see, e.g.,~\cite{BusLop,DriSta94,DenYat,GRS}. Related ideas work for any $r>2$ but becomes very technically demanding, see~\cite{FidLop11}. We briefly remind it here for the reader's convenience.

%

Let $(\mu_1,\mu_2)$ be Nikishin, i.e., 
\begin{equation}\label{eq:2Nikishin}
	d\mu_2(x) = m_{\sigma_2}(x) d\mu_1(x)
\end{equation}
with $\supp\,\mu_1\subseteq\Delta_1$, $\supp\,\sigma_2\subseteq\Delta_2$, $\mathring{\Delta}_1\cap\mathring{\Delta}_2=\varnothing$. 
 
Suppose we are given $(n_1,n_2)\in\bbN^2$ that does not satisfy~\eqref{eq:simpleN}. In fact, we can assume $n_1\le n_2-1$. 
Because of~\eqref{eq:strippingM}, we can rewrite~\eqref{eq:2Nikishin} as
\begin{equation}\label{eq:2NikishinFlipped}
	d\mu_1(x) = (b_1-x-a_1^2 m_{\sigma_2^{(1)}}(x)) d\mu_2(x).
\end{equation}

\begin{lemma}\label{lem:perturbation}
	Let $\bm{\mu}:=(\mu_1,\mu_2)$ and $\widetilde{\bm{\mu}}:=(\widetilde{\mu}_1,{\mu}_2)$, where
	\begin{equation}
		d\widetilde{\mu}_1(x) = d\mu_1(x) + \Big(\sum_{j=0}^s k_j x^j\Big) d\mu_2(x).
	\end{equation}
	Then for any $\bm{n}=(n_1,n_2)\in\bbN^2$ with $n_1 \le n_2-s$, we have $\det\bm{H}_{\bm{n}} = \det \widetilde{\bm{H}}_{\bm{n}}$, so $\bm{n}$ is normal with respect to $\bm{\mu}$ if and only if it is normal with respect to $\widetilde{\bm{\mu}}$. 
\end{lemma}
\begin{proof}
	Note that $H_{n_2,|\bm{n}|}^{(2)}$ block is common for both $\bm{H}_{\bm{n}}$ and $\widetilde{\bm{H}}_{\bm{n}}$. Now note that the $k$-th row of  $\widetilde{H}_{n_1,|\bm{n}|}^{(1)}$ ($k=1,\ldots,n_1$) is 
	\begin{multline}
		\big(\widetilde{c}_{k-1}^{(1)}, \widetilde{c}_{k}^{(1)},\ldots, \widetilde{c}_{k-2+|\bm{n}|}^{(1)}\big)
		=
		\big({c}_{k-1}^{(1)}, {c}_{k}^{(1)},\ldots, {c}_{k-2+|\bm{n}|}^{(1)}\big)
		\\
		+\sum_{j=0}^s k_j
		\big({c}_{k+j-1}^{(2)}, {c}_{k+j}^{(2)},\ldots, {c}_{k+j-2+|\bm{n}|}^{(2)}\big).
	\end{multline}
	But each row of moments in the last sum can be canceled in $\det \widetilde{\bm{H}}_{\bm{n}}$ by subtracting a multiple of the corresponding row from $H_{n_2,|\bm{n}|}^{(2)}$ (use $n_1 \le n_2-s$). This reduces $\det \widetilde{\bm{H}}_{\bm{n}}$  to $\det {\bm{H}}_{\bm{n}}$.
\end{proof}
\begin{remark}
	It is not hard to show that in the setting of Lemma~\ref{lem:perturbation}, $\bm{\mu}$ and $\widetilde{\bm{\mu}}$ share all the type II orthogonal polynomials for indices $(n_1,n_2)$ with $n_1\le n_2-s-1$. 
\end{remark}

\begin{theorem}
	Let $(\mu_1,\mu_2)$ be Nikishin. Then every index $\bm{n}\in\bbN^2$ is normal.
\end{theorem}
\begin{proof}
	For indices with $n_1\ge n_2-1$ this is proven in Theorem~\ref{thm:Nikishin}. For indices with $n_1 \le n_2-1$ (note that indices with $n_1=n_2-1$ fall into both categories), we first observe that $\bm{H}_{n_1,n_2}$ is equal to $\bm{H}_{n_2,n_1}$ for the system $(\mu_2,\mu_1)$. But by~\eqref{eq:2NikishinFlipped}, $\det \bm{H}_{n_2,n_1}$ for $(\mu_2,\mu_1)$ is the same as for  $(\mu_2,\widetilde{\mu}_1)$, where
	$$
	d\widetilde{\mu}_1(x) = -a_1^2 m_{\sigma_2^{(1)}}(x) d\mu_2(x).
	$$
	But this means that $(\mu_2,\widetilde{\mu}_1) = \calN(\mu_2,-a_1^2 \sigma_2^{(1)})$ is Nikishin (observe that $\supp\,\mu_2\subseteq\Delta_1$ and $\supp\,\sigma_2^{(1)}\subseteq \Delta_2$). By Theorem~\ref{thm:Nikishin}, every index $(n_2,n_1)$ with $n_2\ge n_1-1$ is normal for $(\mu_2,\widetilde{\mu}_1)$, which implies that every index with $n_1 \le n_2-1$ is normal for $(\mu_1,\mu_2)$.
\end{proof}

\section{Multiple (Laurent) orthogonal polynomials on the unit circle (MOPUC)}

\subsection{Main definitions}
\hfill\\
 
Let us introduce the notion of multiple Laurent orthogonal polynomials on the unit circle, as introduced in~\cite{KVMLOPUC}. 

For the rest of the paper we choose a branch of the square root function 
\begin{equation}\label{eq:sqrt}
	z^{1/2}=|z|^{1/2} \exp({i \arg_{[t_0,t_0+2\pi)}} z/2)
\end{equation}
for some fixed $t_0\in\bbR$.

Let $r\ge 1$ and consider a system of measures  on the unit circle  $\bm{\mu} = (\mu_1,\dots,\mu_r) \in \mathcal{M}(\bbT)^r$. 

\begin{definition}\label{def:MLOPUC}
We say that an index $\bm{n}\in\bbN^r$ is {$\phi$-}normal with respect to $\bm{\mu}$ if there exists a unique $\phi_{\bm{n}}$ 
 of the form
\begin{equation}\label{eq:monic type II mod}
	\phi_{\bm{n}}(z) = z^{\abs{\bm{n}}/2} + \kappa_{\abs{\bm{n}}/2-1}z^{\abs{\bm{n}}/2-1} + \ldots + \kappa_{-\abs{\bm{n}}/2}z^{-\abs{\bm{n}}/2},
\end{equation} 
that satisfies the orthogonality relations
\begin{equation}\label{eq:type II def mod}
	\int \phi_{\bm{n}}(z)z^{-k}d\mu_j(z) = 0, \quad k = -n_j/2,-n_j/2+1,\ldots,n_j/2-1.
\end{equation}
\end{definition}

It is easy to verify~\cite{KVMLOPUC} that $\phi$-normality of $\bm{n}\ne \bm{0}$ is equivalent to the invertibility of the following $|\bm{n}|\times |\bm{n}|$ matrix 
\begin{equation}
	\label{eq:DOPUC}
	\bm{T}_{\bm{n}} = 
	\begin{pmatrix}
		T^{(1)}_{n_1,|\bm{n}|} \\
		\vdots \\
		T^{(r)}_{n_r,|\bm{n}|}
	\end{pmatrix},
\end{equation}
where
where $T^{(j)}_{n_j,|\bm{n}|}$ is the $n_j\times |\bm{n}|$ matrix given by
\begin{equation}\label{eq:G}
	T^{(j)}_{n_j,|\bm{n}|} =
\begin{pmatrix}
	 c^{(j)}_{-\abs{\bm{n}}/2+n_j/2} &c^{(j)}_{-\abs{\bm{n}}/2+n_j/2+1} & \cdots & c^{(j)}_{\abs{\bm{n}}/2+n_j/2-1} \\
	c^{(j)}_{-\abs{\bm{n}}/2+n_j/2-1} & c^{(j)}_{-\abs{\bm{n}}/2+n_j/2} & \cdots & c^{(j)}_{\abs{\bm{n}}/2+n_j/2-2} \\
	\vdots & \vdots& \ddots & \vdots \\ 
	c^{(j)}_{-\abs{\bm{n}}/2-n_j/2+1} & c^{(j)}_{-\abs{\bm{n}}/2-n_j/2+2} &\cdots & c^{(j)}_{\abs{\bm{n}}/2-n_j/2}
\end{pmatrix}
\end{equation}
with 
$c^{(j)}_k = \int z^k d\mu_j(z)$.

\subsection{AT systems on the unit circle}
\hfill\\

\begin{definition}\label{def:ATOPUC}
	Let  $\Gamma= \{e^{i\theta}:\alpha\le \theta\le \beta\}\subseteq\bbT$ be an arc  with $0<\beta-\alpha\le 2\pi$ and $\mu\in\calM(\Gamma)$ of infinite support. 
	For each $j = 1,\dots,r$, let $d\mu_j(\theta)=w_j(\theta)d\mu(\theta)$ be a positive measure that is absolutely continuous with respect to $\mu$.
	We say that $\bm{\mu}$ is an AT system on $\Gamma$ at $\bm{n}\in\bbN^r$ if 
	\begin{multline}\label{eq:ATOPUC}
		\big\{ z^{-(n_1-1)/2} w_1(z),z^{-(n_1-3)/2} w_1(z),\ldots, z^{(n_1-1)/2} w_1(z),
		\\
		 z^{-(n_2-1)/2} w_2(z),z^{-(n_2-3)/2} w_2(z),\ldots, z^{(n_2-1)/2} w_2(z),
		 \ldots,
		 \\
		  z^{-(n_r-1)/2} w_r(z),z^{-(n_r-3)/2} w_r(z),\ldots, z^{(n_r-1)/2} w_r(z)
		\big\}
	\end{multline}
	is Chebyshev on $\Gamma$.
\end{definition}
Note that $w_j$'s are assumed to be real-valued, but the functions in~\eqref{eq:ATOPUC} are complex-valued. By taking linear combinations one can replace them with trigonometric real-valued functions, see~\cite{KVMLOPUC}. 

\begin{theorem}
	If $\bm{\mu}$ is AT at $\bm{n}$, then  $\bm{n}$ is $\phi$-normal for $\bm{\mu}$.
\end{theorem}
\begin{proof}
	This is proved in~\cite{KVMLOPUC}.
\end{proof}

\subsection{Nikishin systems on the unit circle: same-parity indices with non-increasing components are normal}\label{ss:NikishinOPUC}
\hfill\\

Now let us define Nikishin systems on the unit circle. 
Recall the definition~\ref{eq:caratheodory function} of Carath\'{e}odory functions and that~\eqref{eq:CaraReal} holds for any $\mu\in\mathcal{M}(\bbT)$.

Given two arcs $\Gamma_1$ and $\Gamma_2$ with $\mathring{\Gamma}_1\cap\mathring{\Gamma}_2=\varnothing$, and two measures $\sigma_1\in\calM(\Gamma_1)$,  $\sigma_2\in\calM(\Gamma_2)$, denote $\la  \sigma_1,\sigma_2 \ra$ to be the finite measure 
on $\Gamma_1$ given by
\begin{equation}\label{eq:bracketsOPUC}
	d \la  \sigma_1,\sigma_2 \ra(z) := iF_{\sigma_2}(z) d\sigma_1(z). 
\end{equation}
Just as for the Nikishin system on $\bbR$, here by writing $\la  \sigma_1,\sigma_2 \ra$ we implicitly assume that such a 
finite measure exists. 

\begin{definition}
	We say that $\bm{\mu}=(\mu_1,\ldots,\mu_r)\in\calM(\bbT)^r$ forms a Nikishin system on the unit circle generated by $(\sigma_1,\ldots,\sigma_r)$ (we then write $\bm{\mu}=\calN(\sigma_1,\ldots,\sigma_r)$), if there is
	a collection of arcs $\Gamma_j$, $j=1,\ldots,r$ such that
	\begin{align}
		\label{eq:nonintersectingOPUC}
		& \mathring{\Gamma}_j \cap \mathring{\Gamma}_{j+1} = \varnothing, \qquad j=1,\ldots,r-1, \\
		\label{eq:existst0}
		& \bbT\setminus \bigcup_{j=1}^r \Gamma_j \ne \varnothing,
	\end{align}
	and measures $\sigma_j \in\calM(\Gamma_j)$, $j=1,\ldots,m$, so that
	\begin{equation}\label{eq:NikishinOPUC}
		\mu_1 = \sigma_1,  \mu_{2} = \la \sigma_1, \sigma_{2}\ra, \mu_{3} = \la \sigma_1, \la \sigma_{2},\sigma_3\ra \ra, \ldots, \mu_r = \la \sigma_1, \la \sigma_{2},\la \sigma_3,\ldots,\sigma_r\ra\ra \ra.
	\end{equation}
\end{definition}

We need to choose $t_0$ for the square root branch~\eqref{eq:sqrt}. We do so in such a way that $t_0\notin \Gamma_j$ for any $j$, which is possible from~\eqref{eq:existst0}.


Introduce the notation
\begin{equation}\label{eq:mtkOPUC}
	\begin{aligned}
		m_{w^j\mu} (z) &:=  \int \frac{w^j d\mu(w)}{z-w}, 
		\\
		m_{\la w^j \sigma_1, \sigma_{2}\ra} (z) &: =  i\int_{\Gamma_1} \frac{w^j F_{\sigma_2}(w) d\sigma_1(w)}{z-w}.
	\end{aligned}
\end{equation}

The proof of the following theorem follows the same strategy as our proof of Theorem~\ref{thm:Nikishin} with some minor complications. 
\begin{theorem}\label{thm:NikishinOPUC}
	Let $\bm{\mu}=(\mu_1,\ldots,\mu_r) = \calN(\sigma_1,\ldots,\sigma_r)$ be a Nikishin system on the unit circle. Then
	\begin{enumerate}[(i)]
		\item Let $\bm{n}=(n_1,\ldots,n_r)\in\bbN^r$ with each $n_j$ odd and
\begin{equation}\label{eq:simpleNopuc1}
	n_1\ge n_2 \ge \ldots \ge n_r.
\end{equation}

	Then $\bm{\mu}$ is AT on $\Gamma_1$ at $\bm{n}$. 
	\item Let $\bm{n}=(n_1,\ldots,n_r)\in\bbN^r$ with each $n_j$ even and
\begin{equation}\label{eq:simpleNopuc2}
	n_1\ge n_2 \ge \ldots \ge n_r.
\end{equation}
	Then $\bm{\mu}$ is AT on $\Gamma_1$ at $\bm{n}$. 
	\end{enumerate}
	As a result, any index in {\it (i)} or {\it (ii)} is $\phi$-normal. 
\end{theorem}
\begin{proof}
	 {\it (i)} Each $\mu_j$ is absolutely continuous with respect to $\sigma_1$, with the Radon--Nikodym derivatives 
	 $d\mu_j/d\sigma_1$ equal to 
	 $1$, $iF_{\sigma_2}(x)$, $iF_{\la\sigma_2,\sigma_3\ra}(x)$, $\ldots$, 
	 $iF_{\la\sigma_2,\ldots,\sigma_r\ra }(x)$, respectively. Let $\bm{n}$ satisfies~\eqref{eq:simpleNopuc1}, and each $n_j=2k_j-1$ is odd. 
	 To show the AT property at $\bm{n}$ we need to show nonvanishing of  
	\begin{multline}\label{eq:UOPUC}
		U_{\bm{z}}\big(
		z^{-k_1+1} ,z^{-k_1+2},\ldots, z^{k_1-1},
		z^{-k_2+1} iF_{\sigma_2},z^{-k_2+2} iF_{\sigma_2},\ldots, z^{k_2-1} iF_{\sigma_2},
		\ldots,
		\\
		z^{-k_r+1} F_{\la\sigma_2,\ldots,\sigma_r\ra },\ldots, z^{k_r-1} F_{\la\sigma_2,\ldots,\sigma_r\ra }
		\big),
	\end{multline}
	with an arbitrary $\bm{z}\in\Gamma_1^{|\bm{n}|}$. 
	Let $z_j = e^{i\theta_j}$ with the ordering $t_0 \le \theta_1 <\theta_2 <\ldots < \theta_{|\bm{n}|}<t_0+2\pi$.
	Denote~\eqref{eq:UOPUC} by simply $U_{\bm{z}}$. We will show a slightly stronger statement 
	 		\begin{equation}\label{eq:bigSign}
	 			\frac{U_{\bm{z}}}{|U_{\bm{z}}|}
	 					= i^{l_{\bm{n}}}
	 		\end{equation}
	 		for some integer $l_{\bm{n}}\in\bbZ$ independent of $\bm{z}$, so that, in particular, $U_{\bm{z}}\ne 0$.
	
	We prove the statement using an inductive argument. For $r=1$, we have $n=2k_1-1$ and we easily reduce $U_{\bm{z}}$ to the Vandermonde determinant to get
	\begin{align}
	U_{\bm{z}} &= \prod_{j=1}^{2k_1-1} z_j^{-k_1+1} \prod_{j>k} (z_j-z_k)
	\\
	&=
	\prod_{j=1}^{2k_1-1} z_j^{-k_1+1} i^{n(n-1)/2} \prod_{j=1}^n z_j^{(n-1)/2} \prod_{j>k} |z_j-z_k|
	\\
	&=
	 i^{(2k_1-1)(k_1-1)} 
	 \prod_{j>k} |z_j-z_k|,
	\end{align}
	where in the middle line we used the elementary equality
	\begin{equation}\label{eq:zSin}
		e^{i\theta_j}-e^{i\theta_k} = 2i e^{i(\theta_j+\theta_k)/2} \sin\tfrac{\theta_j-\theta_k}{2}=i z_j^{1/2} z_k^{1/2} |	z_j-z_k |,
	\end{equation}
	because of the ordering $\theta_j>\theta_k$.

	Suppose~\eqref{eq:bigSign} is proven for Nikishin systems on the unit circle consisting of $r-1$ measures for all $\bm{n}$ satisfying \textit{(i)}. 	
	
	First, rewrite $U_{\bm{z}}$ as
	\begin{multline}\label{eq:UOPUCpositive}
		\Big(\prod_{j=1}^{|\bm{n}|} z_j^{-k_1+1} \Big) \, U_{\bm{z}}\big(
		1 ,z, z^2,\ldots, z^{2k_1-2},
		z^{k_1-k_2} F_{\sigma_2},z^{k_1-k_2+1} F_{\sigma_2},\ldots, \\
		z^{k_1+k_2-2} F_{\sigma_2},
		\ldots,
		z^{k_1-k_r} F_{\la\sigma_2,\ldots,\sigma_r\ra },\ldots, z^{k_1+k_r-2} F_{\la\sigma_2,\ldots,\sigma_r\ra }
		\big).
	\end{multline}
	
	By~\eqref{eq:FandM},
	\begin{equation}
		z^j F_\mu(z) = z^{j} - 2 z^{j+1} m_\mu(z)  
	\end{equation}
	for any $j$. Using $\frac{z^{j+1}}{z-w_2} = \sum_{s=0}^{j} z^s w_2^{j-s} + \frac{w_2^{j+1}}{z-w_2}$, we can therefore write 
	\begin{align*}
		z^j F_{\la\sigma_2,\ldots,\sigma_d\ra }(z) 
		&=z^{j} - 2 z^{j+1} m_{\la\sigma_2,\ldots,\sigma_d\ra }(z) 
		\\
		&= z^j - 2i z^{j+1} \int_{\Gamma_2}  \frac{F_{\la\sigma_3,\ldots,\sigma_d\ra}(w_2) }{z-w_2}d\sigma_2(w_2)
		\\
		&= z^j - \sum_{s=0}^{j} v_{j,d,s} z^s 
		-2i \int_{\Gamma_2} \frac{w_2^{j+1} F_{\la\sigma_3,\ldots,\sigma_d\ra}(w_2)  }{z-w_2}d\sigma_2(w_2)
		\\
		&= z^j - \sum_{s=0}^{j} v_{j,d,s} z^s 
		-2  m_{\la w^{j+1} \sigma_2,\ldots,\sigma_d\ra }(z) ,
	\end{align*}
	where 
	$v_{j,d,s} = 2\int_{\Gamma_2} w_2^{j-1-s} d\la\sigma_2,\ldots,\sigma_d\ra(w_2)$ are	
	constants (recall the notation~\eqref{eq:mtkOPUC}).
	
	
	Since $k_1 \ge k_d$ for any $d$, we can perform elementary row operations in the determinant~\eqref{eq:UOPUCpositive} to 	
	reduce it to 
	\begin{multline}\label{eq:U2OPUC}
		(-2)^{|\bm{n}|-n_1} \Big(\prod_{j=1}^{|\bm{n}|} z_j^{-k_1+1} \Big) \, U_{\bm{z}}\big(
		1,z,z^2,\ldots, z^{2k_1-2},
		m_{w^{k_1-k_2+1} \sigma_2}, m_{w^{k_1-k_2+2} \sigma_2},\\
		\ldots, m_{ w^{k_1+k_2-1}\sigma_2},
		\ldots,
		m_{\la w^{k_1-k_r+1}\sigma_2,\ldots,\sigma_r\ra },\ldots,  m_{\la w^{k_1+k_r-1}\sigma_2,\ldots,\sigma_r\ra }
		\big).
	\end{multline}
	Now we apply the generalized Andreief identity, Lemma~\ref{lem:det}{\it (i)}. We treat~\eqref{eq:U2OPUC} as the left-hand side of~\eqref{eq:Andreief} with $N=|\bm{n}|$, $M=|\bm{n}|-n_1$, $N-M = n_1$, $\bm{A} = (z_k^{j})_{0\le j \le 2k_1-2,1\le k \le |\bm{n}|}$, $\mu = \sigma_2$, $g_k(w_2) = \frac{w_2^{k_1}}{z_k - w_2}$. Then the right-hand side of~\eqref{eq:Andreief} becomes
	\begin{multline}\label{eq:inductiveOPUC}
		\tfrac{(-2)^{|\bm{n}|-n_1} \prod_{j=1}^{|\bm{n}|} z_j^{-k_1+1}}{(|\bm{n}|-n_1)!}\int_{\Gamma_2^{|\bm{n}|-n_1}}
		U_{\bm{z}}\Big(1,z,z^2,\ldots, z^{2k_1-2},\tfrac{y_1^{k_1}}{z-y_1},\ldots,\tfrac{y_{|\bm{n}|-n_1}^{k_1}}{z-y_{|\bm{n}|-n_1}}
		\Big)
		\\
		\times 
		U_{\bm{y}}\Big( z^{-k_2+1}, \ldots,z^{k_2-1}, z^{-k_3+1}iF_{\sigma_3}, z^{-k_3+2} iF_{\sigma_3},\ldots, z^{k_3+1} iF_{\sigma_3}, \ldots,
		\\
		 z^{-k_r+1}iF_{\la\sigma_3,\ldots,\sigma_r\ra }, z^{-k_r+2} iF_{\la\sigma_3,\ldots,\sigma_r\ra },
		\ldots,
		 z^{k_r-1} iF_{\la\sigma_3,\ldots,\sigma_r\ra } \Big)
		d^{{|\bm{n}|-n_1}}\sigma_2(\bm{y}).
	\end{multline}
	Observe that the integrand in the last expression is invariant under permutation of $y_j$'s. Therefore we can restrict the domain of integration $\Gamma_2^{|\bm{n}|-n_1}$ to the ordered chamber $\arg y_1<\arg y_2<\ldots<\arg y_{|\bm{n}|-n_1}$ (the argument is taken with the branch cut at $e^{it_0}$ as usual) by simply removing the factor $\tfrac{1}{(|\bm{n}|-n_1)!}$. 
	
	By the induction hypothesis, the system $\mathcal{N}(\sigma_2,\sigma_3,\ldots,\sigma_r)$ is AT at the location $(n_2,n_3,\ldots,n_3)$ (note that this truncated index also satisfies~\eqref{eq:simpleNopuc1}) with the complex phase independent of $\bm{y}$. Therefore we simply need to examine the complex phase of the remaining factors in ~\eqref{eq:inductiveOPUC}.
	
	This can be by reducing the first integrand in ~\eqref{eq:inductiveOPUC} to the Cauchy--Vandermonde determinant (see Lemma~\ref{lem:det}\textit{(ii)}) by extracting $y_j^{k_1}$ from the corresponding rows:
	\begin{align}
		(-2)^{|\bm{n}|-n_1} &\prod_{j=1}^{|\bm{n}|} z_j^{-k_1+1}
		 U_{\bm{z}}\Big(1,z,z^2,\ldots, z^{2k_1-2},\tfrac{y_1^{k_1}}{z-y_1},\ldots,\tfrac{y_{|\bm{n}|-n_1}^{k_1}}{z-y_{|\bm{n}|-n_1}}
		\Big)
		\\
		&=
		(-2)^{|\bm{n}|-n_1} \prod_{j=1}^{|\bm{n}|} z_j^{-k_1+1}
		\prod_{j=1}^{|\bm{n}|-n_1} y_j^{k_1} (-1)^{p_{\bm{n}}}
		\\
		&\qquad 
		\prod_{j<k} (z_k-z_j) \prod_{j<k} (y_k-y_j) \prod_{j,k} (z_j-y_k)^{-1}
		\\
		&=(-2)^{|\bm{n}|-n_1} (-1)^{p_{\bm{n}}} \prod_{j=1}^{|\bm{n}|-n_1} y_j^{k_1}
		\\
		&\qquad \times i^{|\bm{n}|(|\bm{n}|-1)/2} \prod_{j=1}^{|\bm{n}|} {z_j^{(|\bm{n}|-1)/2}}\prod_{j<k} |z_k-z_j|
		\\
		&\qquad \times i^{(|\bm{n}|-n_1)(|\bm{n}|-n_1-1)/2} \prod_{j=1}^{|\bm{n}|-n_1} {y_j^{(|\bm{n}|-n_1-1)/2}}
			 \prod_{j<k} |y_k-y_j|
			 \\
			&\qquad  \times
		 i^{-(|\bm{n}|-n_1)|\bm{n}|} \prod_j {y_j^{-|\bm{n}|/2}} \prod_j z_j^{-(|\bm{n}|-n_1)/2}
		 \prod_{j,k} |z_j-y_k|^{-1},
	\end{align}
	with $p_{\bm{n}}:={n_1(|\bm{n}|-n_1)+\tfrac{(|\bm{n}|-n_1)(|\bm{n}|-n_1-1)}2}$, where we have used~\eqref{eq:zSin} in the last three lines. Note that the factors involving $z_j$'s and $y_j$'s cancel, which proves~\eqref{eq:bigSign} for $r$ measures and completes the induction step.
	
	\textit{(ii)} The argument for the case of all even $n_j=2k_j$ are very similar.
\end{proof}
\begin{remark}
	It is not clear to us if Nikishin system on the unit circle are AT at any index $\bm{n}$ with a mixture of even and odd $n_j$'s, nor whether such an index is necessarily $\phi$-normal.
\end{remark}

\subsection{Nikishin systems on the unit circle, $r=2$ case: all same-parity indices are normal}\label{ss:NikishinOPUCr=2}
\hfill\\

Here we prove that if $r=2$ then the requirement~\eqref{eq:simpleNopuc1}/\eqref{eq:simpleNopuc2} can be removed under a natural requirement on $F_{{\sigma_2}}$. 
It is an interesting open question to show this for $r\ge 3$.

\begin{theorem}
	Let $(\mu_1,\mu_2)=\mathcal{N}(\sigma_1,\sigma_2)$ be Nikishin on $\bbT$. Suppose that $F_{\sigma_2}(z)$ is nonzero on $\Gamma_1$. Then
	\begin{enumerate}[(i)]
		\item  Any $\bm{n}=(n_1,n_2)\in\bbN^2$ with $n_1$ and $n_2$ odd is $\phi$-normal.
		
		\item Any $\bm{n}=(n_1,n_2)\in\bbN^2$ with $n_1$ and $n_2$ even is $\phi$-normal.
	\end{enumerate}
\end{theorem}
\begin{remark}
	If $F_{\sigma_2}(z)$ has a zero on $\Gamma_1$, then $\mu_2 = \la  \sigma_1,\sigma_2 \ra$, see~\eqref{eq:bracketsOPUC}, is typically not a \textit{sign-definite} measure any more. So the condition $F_{\sigma_2}(z)\ne 0$ on $\Gamma_1$ is not only very natural (observe that on the real line, $m_{\sigma_2}(x)\ne 0$ on $\bbR\setminus\Delta_2$ automatically), but it is quite a miracle that the proof of Theorem~\ref{thm:NikishinOPUC} does not require it.
\end{remark}
\begin{remark}
	As we discussed in Section~\ref{ss:OPUC} (see also~\cite{OPUC1} and references therein for more details), $F_{\sigma_2}(z)$ can have at most one zero on $\bbT\setminus\Delta_2$. At that zero $F_{\sigma_2^{(-1)}}$ (see the proof below) will therefore have a pole, which means that $\sigma_2^{(-1)}$ has a point mass there. This strongly suggests that $\phi$-normality of indices with $n_1 < n_2$ need not hold if $F_{\sigma_2}(z)$ vanishes within $\Gamma_1$.
\end{remark}
\begin{proof}
We prove the statement for \textit{(i)} only, as \textit{(ii)} can be proved in the same way. Since $(\mu_1,\mu_2)$ be Nikishin on $\bbT$, we have 
\begin{equation}\label{eq:2NikishinOPUC}
	d\mu_2(z) = iF_{\sigma_2}(z) d\mu_1(z)
\end{equation}
with $\supp\,\mu_1\subseteq\Gamma_1$, $\supp\,\sigma_2\subseteq\Gamma_2$, $\mathring{\Gamma}_1\cap\mathring{\Gamma}_2=\varnothing$. 

Normality of $(n_1,n_2)$ with $n_1\ge n_2$ is proved in Theorem~\ref{thm:NikishinOPUC}. 
Suppose we are given $(n_1,n_2)\in\bbN^2$ with $n_1\le n_2$ (note that diagonal entries fall into both categories). 
Because of~\eqref{eq:CaraInverse}, we can rewrite~\eqref{eq:2NikishinOPUC} as
\begin{equation}\label{eq:2NikishinFlipped}
	d\mu_1(z) = -iF_{\sigma_2^{(-1)}}(z) d\mu_2(z).
\end{equation}
We have $\supp\mu_2\subseteq\Gamma_1$, and since $F_{\sigma_2}\ne 0$ on $\Gamma_1$, we obtain that $\supp\,\sigma_2^{(-1)}$ belongs to $\Gamma_2$ plus potentially one extra point outside of $\Gamma_1\cup\Gamma_2$. In either case we can see that $(\mu_2,\mu_1)$ is a Nikishin system $\mathcal{N}(\mu_2,\sigma_2^{(-1)})$. Applying Theorem to index $(n_2,n_1)$ for this ``flipped'' system~\ref{thm:NikishinOPUC} shows that $(n_1,n_2)$ is $\phi$-normal for $(\mu_1,\mu_2)$.
\end{proof}


\bibsection


\end{document}